\documentclass{amsart}

\def\PSL{\operatorname{PSL}}
\def\PGL{\operatorname{PGL}}
\def\GL{\operatorname{GL}}
\def\AGL{\operatorname{AGL}}

\begin{document}

\title[Automorphisms of Riemann surfaces]
{Automorphism groups of Riemann surfaces of genus $p+1$, where $p$ is prime}

%    Information for first author
\author{Mikhail Belolipetsky}
\address{Mikhail Belolipetsky\newline
Sobolev Institute of Mathematics, Koptyuga 4, 630090 Novosibirsk, Russia\newline
Einstein Institute of Mathematics, Hebrew University, Jerusalem 91904, Israel}
\email{mbel@math.huji.ac.il}
%    Information for second author
\author{Gareth A.~Jones}
\address{Gareth A.~Jones\newline
School of Mathematics, University of Southampton, Southampton SO17 1BJ,\newline
United Kingdom}
\email{G.A.Jones@maths.soton.ac.uk}


\begin{abstract}
We show that if $\mathcal S$ is a compact Riemann surface of genus $g=p+1$, where
$p$ is prime, with a group of automorphisms $G$ such that $|G|\geq\lambda(g-1)$
for some real number $\lambda>6$, then for all sufficiently large $p$
(depending on $\lambda$), $\mathcal S$ and $G$ lie in one of six infinite sequences
of examples. In particular, if $\lambda=8$ then this holds for all $p\geq 17$.
\end{abstract}

\maketitle

\section{\bf Introduction}

A compact Riemann surface $\mathcal S$ of genus $g\geq 2$ has at most $84(g-1)$
automorphisms, and among the most interesting surfaces are those with a group
$G$ of automorphisms which is relatively large compared with $g$. The general
problem of determining all such surfaces $\mathcal S$ and groups $G$ is very
difficult, but it tends to be easier when the Euler characteristic
$\chi=2(1-g)$ of $\mathcal S$ has a simple numerical form. Here we will consider
the simplest case, where $g=p+1$ for some prime $p$. In this situation, Accola
[{\bf 2}] has determined the possibilities for $G$ where $|G|\geq 8(g+1)$, and
we will extend his results to the case $|G|\geq \lambda(g-1)$ for each
$\lambda>6$. We will show that if $p$ is sufficiently large (depending on
$\lambda$), then $\mathcal S$ and $G$ lie in one of six easily described infinite
families. Our method is a combination of the techniques used by Accola and
those developed by us in [{\bf 3}], where bounds were obtained for automorphism
groups of compact arithmetic Riemann surfaces. Our main result is the following
(where $\mathcal H$ denotes the hyperbolic plane and $\Gamma(l,m,n)$ denotes a
triangle group):

\medskip

T{\small HEOREM} 1. {\it For each real number $\lambda>6$ there is a constant
$c_{\lambda}$ with the following property. Let ${\mathcal S}$ be a compact Riemann
surface of genus $g=p+1$ for some prime $p\geq c_{\lambda}$, and suppose that
some subgroup $G\leq{\rm Aut}\,({\mathcal S})$ has order $|G|\geq \lambda(g-1)$.
Then
\begin{itemize}
\item[\rm(a)] ${\mathcal S}\cong{\mathcal H}/K$ and $G\cong\Gamma/K$ for some
Fuchsian group $\Gamma$ and normal surface subgroup $K$ of $\Gamma$, where one
of the following holds:
\vskip4pt
\item[]{\rm(i)} $|G|=12(g-1)$ where $p\equiv 1$ {\rm mod}~$(3)$,
$\Gamma=\Gamma(2,6,6)$, and $G$ is a split extension of $C_p$ by $C_6\times
C_2$;
\vskip4pt
\item[]{\rm(ii)} $|G|=10(g-1)$ where $p\equiv 1$ {\rm mod}~$(5)$,
$\Gamma=\Gamma(2,5,10)$, and $G$ is a split extension of $C_p$ by $C_{10}$;
\vskip4pt
\item[]{\rm(iii)} $|G|=8(g-1)$ where $p\equiv 1$ {\rm mod}~$(8)$,
$\Gamma=\Gamma(2,8,8)$, and $G$ is a split extension of $C_p$ by $C_8$;
\vskip4pt
\item[]{\rm(iv)} $|G|=8(g+3)$ where $p\equiv 2$ {\rm mod}~$(3)$,
$\Gamma=\Gamma(2,4,n)$ with $n=p+4$, and $G$ is an extension of $C_{n/3}$ by
$S_4$; \vskip4pt
\item[]{\rm(v)} $|G|=8(g+1)$ where $\Gamma=\Gamma(2,4,n)$ with $n=2p+4$,
and $G$ is an extension of $C_{n/2}$ by $D_4$; \vskip4pt
\item[]{\rm(vi)} $|G|=8g$ where $\Gamma=\Gamma(2,4,n)$ with $n=4p+4$,
and $G$ is an extension of $C_n$ by $C_2$. \vskip4pt
\item[\rm(b)] $G={\rm Aut}\,({\mathcal S})$, and $\Gamma$ is the normalizer
$N(K)$ of $K$ in $\PSL(2,{\bf R})$. \vskip4pt
\item[\rm(c)] In cases~(i) and (iii), for each prime $p$ satisfying the
given congruence there are two non-isomorphic surfaces $\mathcal S$, forming a
chiral pair; in case~(ii) there are four surfaces for each $p$, forming two
chiral pairs, and in cases~(iv), (v) and (vi) $\mathcal S$ is unique. In each case,
the group $G$ is determined uniquely (up to isomorphism) by $p$.
\end{itemize}
}
\medskip

We will give more detailed descriptions of these surfaces $\mathcal S$ and groups
$G$ in Section~3. For instance, the surfaces in cases~(iv) and (v) are the
well-known examples constructed by Accola [{\bf 1}] and Maclachlan [{\bf 11}],
while the group $G$ in (vi) is the semidihedral group $SD_n$ of order $2n$.

\medskip

As a specific example of Theorem~1, we will show in Section~6 that one can take
$c_8=17$, so that if $|G|\geq 8(g-1)$ then the conclusions of Theorem~1 are
valid for all $p\geq 17$. This is an extension of Accola's result in [{\bf 2},
chapter 7], where he considered the case $|G|\geq 8(g+1)$ (his `big groups')
and showed that if $p\geq 89$ then only cases~(i), (ii), (iv) and (v) occur.

\medskip

These lower bounds on $p$ are necessary to avoid sporadic examples for low
genera which do not belong to the six infinite families in Theorem~1(a). For
instance, the groups in these families are all solvable, whereas there are also
non-solvable examples for small $p$: these include $G=S_5$ of order $40(g-1)$
for $p=3$, and $G=\PSL_2(13)$ of order $84(g-1)$ for $p=13$, arising from
$\Gamma=\Gamma(2,4,5)$ and $\Gamma(2,3,7)$ respectively.

\medskip

By Dirichlet's Theorem, each of the congruences in Theorem~1(a) is satisfied by
infinitely many primes $p$, so each of cases (i) to (vi) corresponds to an
infinite set of genera $g$. Indeed, we will show in Section~3 that for every
prime $p$ satisfying the given congruences (and not just for all sufficiently
large $p$), there exist surfaces $\mathcal S$ and groups $G$ satisfying the
conditions of cases (i) to (vi). From this, and from the case $\lambda=8$ of
Theorem~1 considered in Section~6, we deduce the following result, which
extends a result of Accola for $p\geq 89$ [{\bf 2}, theorem 7.11]. For each
$g\geq 2$, let $N(g)$ denote the maximum number of automorphisms of a Riemann
surface of genus~$g$.

\medskip

T{\small HEOREM} 2. {\it Let $g=p+1$ for some prime $p\geq 17$.
\begin{itemize}
\item[\rm(a)] If $p\equiv 1$ {\rm mod}~$(3)$ then $N(g) = 12(g-1)$.
\vskip2pt
\item[\rm(b)] If $p\equiv 11$ {\rm mod}~$(15)$ then $N(g) = 10(g-1)$.
\vskip2pt
\item[\rm(c)] If $p\equiv 2, 8$ or $14$ {\rm mod}~$(15)$ then $N(g) = 8(g+3)$.
\end{itemize}}
\medskip

Together with the known results on $N(g)$ for small $g$ (see [{\bf 2}]), namely
$N(14) = 1092$, $N(12) = 120$, $N(8) = 336$, $N(6) = 150$, $N(4) = 120$ and
$N(3) = 168$, this Theorem gives the value of $N(p+1)$ for every prime $p$.

\section{\bf Preliminaries}

In this section we recall some basic facts about Riemann surfaces and their
groups of automorphisms. For more information see [{\bf 5}], [{\bf 8}].

\medskip

By the uniformization theorem [{\bf 5}, chapter IV], each compact Riemann
surface ${\mathcal S}$ of genus $g\geq 2$ is isomorphic to ${\mathcal H}/K$, where
$\mathcal H$ is the hyperbolic plane and $K$ is a torsion-free discrete subgroup of
the group ${\rm Isom}^+({\mathcal H}) = \PSL(2,{\bf R})$ of orientation-preserving
isometries of $\mathcal H$. This group $K$, called the {\it surface group\/}
corresponding to ${\mathcal S}$, is unique up to conjugacy in $\PSL(2,{\bf R})$.

\medskip

Discrete subgroups of $\PSL(2,{\bf R})$ are called {\it Fuchsian groups}. Each
cocompact Fuchsian group $\Gamma$ has a presentation
$$\Gamma=\langle\,\alpha_1,\beta_1,\dots,\alpha_g,\beta_g,\gamma_1,\dots,
\gamma_k\mid\prod_{i=1}^g[\alpha_i,\beta_i]\prod_{j=1}^k\gamma_j=1,\,
\gamma_j^{m_j}=1\,\rangle,$$ with genus $g\geq 0$ and elliptic periods $m_i\geq
2$. We write $\Gamma=\Gamma(g;m_1,\ldots, m_k)$, and we call $(g;m_1,\ldots,
m_k)$ the {\it signature\/} $\sigma$ of $\Gamma$, usually abbreviated to
$(m_1,\ldots, m_k)$ if $g=0$. The order of the elliptic periods is irrelevant,
so we may take $m_1\leq\ldots\leq m_k$. A {\it surface group\/} has signature
$(g;-)$, with $k=0$, and a group with signature $(m_1,m_2,m_3)$ is called a
{\it triangle group}.

\medskip

We define $\mu(\Gamma)$ to be the hyperbolic measure of ${\mathcal H}/\Gamma$, or
equivalently of a fundamental region for $\Gamma$. For cocompact groups this is
finite, and can be expressed in terms of the signature:
$$
\mu(\Gamma)=\mu(g;m_1,\dots ,m_k)=2\pi\left(2g-2+\sum_{j=1}^k \left(1-{1\over
m_j}\right)\right).
$$
The Riemann--Hurwitz formula states that if $\Delta$ is a subgroup of finite
index in $\Gamma$ then $\mu(\Delta)=|\Gamma:\Delta|\cdot\mu(\Gamma)$.

\medskip

The automorphisms of a Riemann surface ${\mathcal S}$ lift to the isometries of
$\mathcal H$ normalizing the surface group $K$, so the automorphism group ${\rm
Aut}\,({\mathcal S})$ is isomorphic to $N(K)/K$ where $N(K)$ is the normalizer of
$K$ in $\PSL(2,{\bf R})$.

\medskip

A {\it surface-kernel epimorphism (SKE)} is an epimorphism $\theta: \Gamma\to
G$ such that $K = {\rm ker}(\theta)$ is a surface group. When $\Gamma$ is
cocompact and  $G$ is finite, this is equivalent to the condition that
$\gamma_j\theta$ has order $m_j$ for $j=1,\ldots, k$. In this situation, the
action of $\Gamma$ on $\mathcal H$ induces a faithful action of $\Gamma/K$ on
${\mathcal H}/K$, so the Riemann surface ${\mathcal S}$ uniformized by $K$ has a group
of automorphisms isomorphic to $G$.

\bigskip

\section{\bf Existence results}

In this section we prove the existence of Riemann surfaces $\mathcal S$ and groups
$G$ described in cases~(i) to (vi) of Theorem~1(a) for all primes $p$
satisfying the given congruences, not just for all sufficiently large $p$. We
will first describe the general method, and then outline the details for the
individual cases.

\medskip

In each case, we will define a SKE $\theta:\Gamma\to G$, so that the surface
group $K={\rm ker}(\theta)$ uniformises a compact Riemann surface ${\mathcal
S}={\mathcal H}/K$. The Riemann--Hurwitz formula, applied to the inclusion
$K\leq\Gamma$ of index $|G|$, gives the genus $g$ of $\mathcal S$, and the fact
that $K$ is normal in $\Gamma$ implies that $G\cong\Gamma/K\leq N(K)/K\cong{\rm
Aut}({\mathcal S})$.

\medskip

For each $K$, $N(K)$ is a Fuchsian group containing $\Gamma$. Now Singerman
[{\bf 14}] has shown that the triangle groups $\Gamma$ in cases (ii), (iv), (v)
and (vi) are maximal Fuchsian groups, so $N(K)=\Gamma$ and hence $G={\rm
Aut}({\mathcal S})$. In cases~(i) and (iii), $\Gamma=\Gamma(2,6,6)$ and
$\Gamma(2,8,8)$ have index $2$ in $N(\Gamma)=\Gamma(2,4,6)$ and
$\Gamma(2,4,8)$, which are maximal Fuchsian groups; we will show that
$N(\Gamma)$ does not normalise $K$ in these cases, so again $N(K)=\Gamma$ and
$G={\rm Aut}({\mathcal S})$.

\medskip

Two SKEs $\Gamma\to G$ have the same kernel if and only if they differ by an
automorphism of $G$, so the number of normal surface subgroups $K$ of $\Gamma$
with quotient group $G$ is equal to the number of orbits of ${\rm Aut}(G)$ on
SKEs $\Gamma\to G$. Only the identity automorphism can fix a SKE, so this
action is semi-regular, and the number of orbits is the number of SKEs divided
by $|{\rm Aut}(G)|$. We can count SKEs $\Gamma\to G$ by counting appropriate
generating sets for $G$, so we can calculate the number of subgroups $K$.

\medskip

Two surface groups uniformise isomorphic Riemann surfaces if and only if they
are conjugate in $\PSL(2,{\bf R})$. Since triangle groups with a given signature
are all conjugate, we can count isomorphism classes of Riemann surfaces $\mathcal
S$ in each case by considering a fixed triangle group $\Gamma$, and counting
conjugacy classes in $\PSL(2,{\bf R})$ of its surface subgroups $K$. If two such
subgroups $K_1$ and $K_2$ satisfy $K_1^{\gamma}=K_2$ for some $\gamma\in
\PSL(2,{\bf R})$ then $N(K_1)^{\gamma}=N(K_2)$; we have seen that
$N(K_1)=N(K_2)=\Gamma$, so $\gamma\in N(\Gamma)$. Thus $K_1$ and $K_2$ are
conjugate in $\PSL(2,{\bf R})$ if and only if they are conjugate in $N(\Gamma)$.
In this action of $N(\Gamma)$ by conjugation on these surface groups, the
stabiliser of each $K$ is $N(K)=\Gamma$, so $K$ lies in an orbit of length
$|N(\Gamma):\Gamma|$. In cases~(i) and (iii), $|N(\Gamma):\Gamma|=2$, so the
subgroups $K$ are conjugate in pairs, with distinct pairs uniformising
non-isomorphic surfaces; in cases~(ii), (iv), (v) and (vi), however,
$N(\Gamma)=\Gamma$, so the subgroups are mutually non-conjugate and their
surfaces are non-isomorphic.

\medskip

\noindent {\it Example\/} (i). Let $\Gamma=\Gamma(2,6,6)$, let $p$ be any prime
such that $p\equiv 1$ mod~$(3)$, and let
$$G=\langle x,y\mid x^p=y^6=z^2=1,\, x^y=x^u,\, [x,z]=[y,z]=1\rangle,$$
where $u$ is a primitive $6$-th root of unity mod~$(p)$. (Since $p\equiv 1$
mod~$(6)$ there are $\phi(6)=2$ mutually inverse choices for $u$, but the
resulting groups $G$ are isomorphic under the mapping $x\mapsto x,\, y\mapsto
y^{-1},\, z\mapsto z$.) Then $G$ is a split extension of a normal subgroup
$P=\langle x\rangle\cong C_p$ by $Q=\langle y,z\rangle\cong C_6\times C_2$. We
also have
$$G=G_1\times Z=G_2\times Z$$
where $Z=\langle z\rangle\cong C_2$ is the centre of $G$, while $G_1=\langle
x,y\rangle$ and $G_2=\langle x, yz\rangle$ are isomorphic subgroups of order
$6p$ and index $2$. It follows that $|{\rm Aut}\,(G)|=2p(p-1)$: automorphisms
must fix $Z$, and can preserve or transpose $G_1$ and $G_2$, so ${\rm
Aut}\,(G)$ has a subgroup ${\rm Aut}\,(G_1)$ of index $2$ under which there are
$p-1$ choices $x^i\;(i\neq 0)$ for the image of $x$, and $p$ choices $x^iy$ for
the image of $y$.

\medskip

Surface-kernel epimorphisms $\theta:\Gamma\to G$ correspond to generating
triples $a,b,c=\gamma_i\theta\;(i=1,2,3)$ of orders $2,6$ and $6$ satisfying
$abc=1$. Each element of $G$ has the unique form $x^iy^jz^k$ where $i\in{\bf
Z}_p,\, j\in{\bf Z}_5$ and $k\in{\bf Z}_2$. The elements of order $2$ are those
of the form $x^iy^3z^k$, together with $z$ (which cannot be a member of such a
triple); the elements of order $6$ are those of the form $x^iy^{\pm 1}z^k$ or
$x^iy^{\pm 2}z$. A little calculation (which we shall omit) then shows that
${\rm Aut}\,(G)$ has four orbits on the required triples $a,b,c$:
\begin{align*}
a=x^iy^3z^k,\quad b&=x^{i'}yz^{1-k},\qquad c=x^{i''}y^2z,\cr
a=x^iy^3z^k,\quad b&=x^{i'}y^{-1}z^{1-k},\quad c=x^{i''}y^{-2}z,\cr
a=x^iy^3z^k,\quad b&=x^{i'}y^2z,\qquad\quad\, c=x^{i''}yz^{1-k},\cr
a=x^iy^3z^k,\quad b&=x^{i'}y^{-2}z,\qquad\;\; c=x^{i''}y^{-1}z^{1-k}.
\end{align*}
In each case there are $p$ choices for $i$ and two choices for $k$, then $p-1$
choices for $i'$ (excluding one value which gives $a=b^3$), and then $i''$ is
uniquely determined by the equation $abc=1$. It follows that there are four
normal surface subgroups $K$ in $\Gamma$ with $\Gamma/K\cong G$, so that the
quotient surfaces ${\mathcal S}={\mathcal H}/K$ have $G\leq{\rm Aut}\,({\mathcal S})$. In
each case, $\mathcal S$ has genus $g=p+1$ since $\mu(\Gamma)=\pi/3$ and $|G|=12p$.
We can label these subgroups and surfaces $K_j$ and ${\mathcal S}_j\;(j=\pm 1, \pm
2\in{\bf Z}_6)$ respectively, as $\gamma_2$ has eigenvalue $u^j$ on $P$
(regarded as a $1$-dimensional vector space over ${\bf Z}_p$). Since $\gamma_2$
is conjugate in $N(\Gamma)=\Gamma(2,4,6)$ to $\gamma_3$, which has eigenvalue
$u^{3-j}$ on $P$, it follows that each $K_j$ is conjugate to $K_{3-j}$, so we
obtain two non-isomorphic surfaces ${\mathcal S}_1\cong{\mathcal S}_2$ and ${\mathcal
S}_{-1}\cong{\mathcal S}_{-2}$. In the extended (orientation-reversing) triangle
group, which contains $\Gamma$ with index $2$, $\gamma_2$ and $\gamma_3$ are
conjugate to their inverses and so each $K_j$ is conjugate to $K_{-j}$; thus
the two surfaces ${\mathcal S}_1$ and ${\mathcal S}_{-1}$ form a chiral pair,
corresponding to complex conjugate algebraic curves.

\medskip

\noindent {\it Example\/} (ii). Here $\Gamma=\Gamma(2,5,10)$ and
$$G=\langle x,y\mid x^p=y^{10}=1,\, x^y=x^u\rangle$$
where $u$ is a primitive $10$-th root of unity mod~$(p)$ for some prime
$p\equiv 1$ mod~$(5)$. Thus $G$ is a split extension of a normal subgroup
$P=\langle x\rangle\cong C_p$ by $Q=\langle y\rangle\cong C_{10}$. Each element
of $G$ has the unique form $x^iy^j$ where $i\in{\bf Z}_p$ and $j\in{\bf
Z}_{10}$; it has order $1$ or $p$ if $j=0$, and otherwise its order is $10/{\rm
hcf}(10,j)$. It follows that there are $4p(p-1)$ SKEs $\theta:\Gamma\to G$:
there are $p$ choices of an element $a=\gamma_1\theta=x^iy^5$ of order $2$, and
then $4(p-1)$ choices of an element $b=\gamma_2\theta=x^{i'}y^{2j}\;(j=\pm
1,\pm 2\in{\bf Z}_5)$ of order $5$, in each case excluding one value of $i'$
for which $\langle a,b\rangle\cong C_{10}$. Now $G$ has $p(p-1)$ automorphisms,
by the argument applied to $G_1$ in Example~(i), so there are
$4p(p-1)/p(p-1)=4$ normal surface subgroups $K$ of $\Gamma$ with $\Gamma/K\cong
G$. The four surfaces ${\mathcal S}={\mathcal H}/K$ have genus $g=p+1$ (since
$\mu(\Gamma)=2\pi/5$ and $|G|=10p$) and $G\leq{\rm Aut}\,({\mathcal S})$; they are
mutually non-isomorphic since $N(\Gamma)=\Gamma$. We can write $K=K_j$ and
${\mathcal S}={\mathcal S}_j\;(j=\pm 1, \pm 2\in{\bf Z}_5)$, where $u^{2j}$ is the
eigenvalue of $\gamma_2$ on $P$. Since $\gamma_2$ is conjugate to its inverse
in the extended triangle group, the surfaces ${\mathcal S}_1$ and ${\mathcal S}_{-1}$
form a chiral pair, as do ${\mathcal S}_2$ and ${\mathcal S}_{-2}$. In fact, more
general results of Streit and Wolfart show that these four surfaces, defined
over the field ${\bf Q}(e^{2\pi i/5})$, are conjugate under the Galois group of
that field (they are the surfaces $X_{n,t,t}$ of [{\bf 15}, theorem~3], with
$q=5$).

\medskip

\noindent {\it Example\/} (iii). Here we imitate Example~(ii), with $p\equiv 1$
mod~$(8)$, $\Gamma=\Gamma(2,8,8)$ and
$$G=\langle x,y\mid x^p=y^8=1, x^y=x^u\rangle$$
where $u$ is a primitive $8$-th root of unity mod~$(p)$. The same arguments as
before show that there are four surface groups $K=K_j$ in $\Gamma$ with
$\Gamma/K\cong G$, distinguished by the eigenvalue $u^j\;(j=\pm 1,\pm 3)$ of
$\gamma_2$ on $P=\langle x\rangle$. Since $\mu(\Gamma)=\pi/2$ and $|G|=8p$ the
four corresponding surfaces ${\mathcal S}={\mathcal S}_j$ all have genus $p+1$, and
since $|N(\Gamma):\Gamma|=2$ they are isomorphic in pairs. Now $\gamma_2$ is
conjugate in $N(\Gamma)=\Gamma(2,4,8)$ to $\gamma_3=(\gamma_1\gamma_2)^{-1}$,
which has eigenvalue $u^{4-j}$, so ${\mathcal S}_1\cong{\mathcal S}_3$ and ${\mathcal
S}_{-1}\cong{\mathcal S}_{-3}$; as before, these two surfaces ${\mathcal S}_1$ and
${\mathcal S}_{-1}$ form a chiral pair.

\medskip

Each of the groups $G$ in Examples~(i), (ii) and (iii) has a normal Sylow
$p$-subgroup $P\cong C_p$, with quotient group $Q=G/P$ isomorphic to $C_6\times
C_2$, $C_{10}$ or $C_8$ respectively. The inverse image $\Delta$ of $P$ in
$\Gamma$ is a surface group of genus $2$, so $\mathcal S$ is a $p$-sheeted
unbranched covering of a Riemann surface ${\mathcal T}={\mathcal H}/\Delta$ of genus
$2$. In each case there is a unique isomorphism class of Riemann surfaces $\mathcal
T$ of genus $2$ with $Q\leq{\rm Aut}\,({\mathcal T})$: in Example~(i) it
corresponds to the algebraic curve $w^2=z^6-1$, with generators of the direct
factors $C_6$ and $C_2$ of $Q$ acting on $\mathcal T$ by $(z,w)\mapsto (e^{\pi
i/3}z,w)$ and $(z,w)\mapsto (z,-w)$ (the hyperelliptic involution); in
Example~(ii) we have $w^2=z^5-1$ with a generator of $Q\cong C_{10}$ acting by
$(z,w)\mapsto (e^{2\pi i/5}z,-w)$; in Example~(iii) we have $w^2=z^5-z$ with a
generator of $Q\cong C_8$ acting by $(z,w)\mapsto (iz,e^{\pi i/4}w)$.

\medskip

\noindent{\it Example\/} (iv). As shown by Accola [{\bf 1}] and Maclachlan
[{\bf 11}], if $n=3m$ for some integer $m\geq 1$ then there a surface-kernel
epimorphism $\theta:\gamma_i\mapsto a,b,c$ from $\Gamma=\Gamma(2,4,n)$ to a
group
$$G=\langle a,b,c\mid a^2=b^4=c^n=abc=1, (c^3)^a=c^{-3}\rangle$$
which is an extension of a normal subgroup $\langle c^3\rangle\cong C_m$ by
$\Gamma(2,4,3)\cong S_4$ (also see [{\bf 4}] for a construction of this group).
Thus $|G|=24m=8n$, and since $\mu(\Gamma)=\pi(n-4)/2n$ it follows that the
surface subgroup $K={\rm ker}(\theta)$ has genus $g=n-3=3(m-1)$, so $G$ is a
group of $8(g+3)$ automorphisms of the surface ${\mathcal S}={\mathcal H}/K$. This
example exists for every genus $g$ divisible by $3$; this includes $g=p+1$ for
any prime $p\equiv 2$ mod~$(3)$, with $n=g+3=p+4$.

\medskip

This situation, together with that described in Example~(v), is part of a wider
investigation by Jones and Surowski [{\bf 9}] of cyclic coverings of the
Platonic maps, or equivalently, cyclic extensions of finite rotation groups
arising as quotients of triangle groups. It follows from their results that if
$n$ is odd (as it is here when $g=p+1$ for odd primes $p\equiv 2$ mod~$(3)$),
there is a unique normal surface subgroup $K$ in $\Gamma=\Gamma(2,4,n)$ with
$G=\Gamma/K$ an extension of $C_m$ by $S_4$ (for even $n=3m$ there may be
several such $K$ and $G$, obtained by replacing the last relation in $G$ with
$(c^3)^a=c^{3u}$ where $u^2\equiv 1$ and $4(1+u)\equiv 0$ mod~$(m)$). The
associated surface $\mathcal S$ is an $m$-sheeted regular cyclic covering of the
sphere, branched over the eight vertices of a cube, with monodromy permutations
$\pm 1$ (in the additive group ${\bf Z}_m$) at alternate vertices.

\medskip

\noindent {\it Example\/} (v). The construction here, also due to Accola and
Maclachlan, is similar to that in Example~(iv). We have $\Gamma=\Gamma(2,4,n)$
where $n=2m$, and
$$G=\langle a,b,c\mid a^2=b^4=c^n=abc=1,\, (c^2)^a=c^{-2}\rangle,$$
an extension of $\langle c^2\rangle\cong C_m$ by $\Gamma(2,4,2)\cong D_4$, of
order $8m=4n$. The surface-kernel $K$ has genus $g=m-1$, so $|G|=8(g+1)$.
Again, it is shown in [{\bf 9}] that $K$ is unique. The surface $\mathcal S$ is an
$m$-sheeted regular cyclic covering of the sphere, branched at $\pm 1$ (with
monodromy permutation $1\in{\bf Z}_m$) and at $\pm i$ (with monodromy
permutation $-1$). This applies to every genus $g\geq 0$, including those of
the form $g=p+1$.

\medskip

\noindent {\it Example\/} (vi). Here $\Gamma=(2,4,n)$ where $n=2m$ for some
even $m$, and
$$G=\langle a,b,c\mid a^2=b^4=c^n=abc=1,\, c^a=c^{m-1}\rangle.$$
Eliminating $b$ we see that
$$G=\langle a,c\mid a^2=c^n=1,\, c^a=c^{m-1}\rangle$$
is the semidihedral group $SD_n$ of order $2n$, an extension of $\langle
c\rangle\cong C_n$ by $\langle a\rangle\cong C_2$ (see Section~5 for the
details). Since $\mu(\Gamma)=\pi(n-4)/2n$ and $|G|=2n$, the corresponding
surface $\mathcal S$ has genus $g=n/4$, so $|G|=8g$. This example is valid for all
$g\geq 2$.

\medskip

We can now make a first step towards proving Theorem~2:

\medskip

C{\small OROLLARY} 3. {\it Let $g=p+1$ for some prime $p$.
\begin{itemize}
\item[\rm(a)] If $p\equiv 1$ {\rm mod}~$(3)$ then $N(g)\geq 12(g-1)$.
\vskip2pt
\item[\rm(b)] If $p\equiv 11$ {\rm mod}~$(15)$ then $N(g)\geq 10(g-1)$.
\vskip2pt
\item[\rm(c)] If $p\equiv 2, 8$ or $14$ {\rm mod}~$(15)$ then $N(g)\geq
8(g+3)$. \vskip2pt \noindent In particular, $N(g)\geq 8(g+3)$.
\end{itemize}
}
\medskip

{\it Proof.} Examples~(i), (ii) and (iv) imply (a), (b) and (c), which cover
all the primes except $p=3$ and $5$. For the last assertion, in~(a) we have
$N(g)\geq 12(g-1)\geq 8(g+3)$ for all $g\geq 9$, and in the remaining case
$g=8$, a surface-kernel epimorphism $\Gamma(2,3,8)\to \PGL_2(7)$ gives $N(g)\geq
336\geq 8(g+3)$. In~(b) we have $N(g)\geq 10(g-1)\geq 8(g+3)$ for all $g\geq
17$, and in the remaining case $g=12$ we can use Example~(iv) to see that
$N(g)\geq 8(g+3)$. Example~(iv) also deals with $g=6$, and for $g=4$ we can use
Bring's curve [{\bf 12}], corresponding to a surface-kernel epimorphism from
$\Gamma(2,4,5)$ to $S_5$, to see that $N(g)\geq 120\geq 8(g+3)$.

\section{\bf Proof of Theorem~1(a)}

Suppose that $G$ and $\mathcal S$ are as in Theorem~1, with $|G|\geq\lambda(g-1)$
for some $\lambda>6$. Since $g\geq 2$ we have $G\cong\Gamma/K$ for some normal
surface subgroup $K$ of a finitely generated Fuchsian group $\Gamma$, with
$$4\pi(g-1)=\mu(K)=|G|\mu(\Gamma)\geq\lambda(g-1)\mu(\Gamma),\eqno(4{\cdot}1)$$
so $\mu(\Gamma)\leq 4\pi/\lambda<2\pi/3$. If $\Gamma$ has signature
$\sigma=(g;m_1,\ldots,m_k)$ then since
$$\mu(\Gamma)=2\pi\left(2g-2+\sum_{j=1}^k \left(1-{1\over
m_j}\right)\right)>0\eqno(4{\cdot}2)$$ it follows by case-by-case analysis (see
[{\bf 2}, \S 7.6] for the case $\mu(\Gamma)<\pi/2$) that $g=0$ and that apart
from a finite set of signatures, such as $(2,6,n)$ for $6\leq n\leq\lfloor
3\lambda/(\lambda-6)\rfloor$, the possible signatures $\sigma$ all lie in three
infinite sequences, namely $(2,3,n)$ for $n\geq 7$, and $(2,4,n)$ and $(2,5,n)$
for $n\geq 5$. (If we had taken $\lambda=6$, we would also have to include
$\sigma=(2,6,n)$ for all $n\geq 6$ and $(3,3,n)$ for all $n\geq 4$; other
infinite sequences of signatures appear for smaller values of $\lambda$.)

\medskip

Let $\Sigma_{\lambda}$ be the finite set consisting of all these signatures
$\sigma$ except $(2,3,n)$ for $n\geq 78$, $(2,4,n)$ for $n\geq 37$, and
$(2,5,n)$ for $n\geq 71$. For $\sigma\in\Sigma_{\lambda}$ we will follow the
method used in [{\bf 3}] for arithmetic Riemann surfaces (though the numerical
details are somewhat different), and then we will use separate arguments for
the remaining three infinite sequences.

\medskip

{\it Case A\/}: $\sigma\in\Sigma_{\lambda}$. By ($4{\cdot}2$), the number
$q=\mu(\Gamma)/4\pi$ is rational and depends only on the signature $\sigma$ of
$\Gamma$, so let us write $q=r/s$ where $r$ and $s$ are coprime positive
integers. Let $p_{\lambda}$ be the largest prime dividing any $r$ for
$\sigma\in\Sigma_{\lambda}$. Now ($4{\cdot}1$) gives $|G|=(g-1)/q=ps/r$; since
this is an integer we have $r=1$ and $|G|=ps$ for all $p>p_{\lambda}$.

\medskip

Let $\Sigma'_{\lambda}=\{\sigma\in\Sigma_{\lambda}\mid r=1\}$ and let
$s_{\lambda}$ be the maximum value of $s$ for $\sigma\in\Sigma'_{\lambda}$. For
$p>s_{\lambda}+1$ we therefore have $|G|=ps$ with $s$ coprime to $p$ and less
than $p+1$, so Sylow's Theorems imply that $G$ has a normal Sylow $p$-subgroup
$P\cong C_p$. The inverse image of $P$ in $\Gamma$ is a normal subgroup
$\Delta$ of $\Gamma$ with quotient $Q=\Gamma/\Delta\cong G/P$ of order $s$.
Since $p$ and $s$ are coprime, the Schur--Zassenhaus Theorem implies that $G$
is a split extension of $P$ by $Q$. Now let $p'_{\lambda}$ be the largest prime
dividing any of the elliptic periods $m_i$ of the signatures
$\sigma\in\Sigma'_{\lambda}$. Then for $p>p'_{\lambda}$, the inclusion
$K\leq\Delta$ induces a smooth $p$-sheeted covering ${\mathcal S}\to{\mathcal T}={\mathcal
H}/\Delta$ of surfaces, so $\Delta$ is a surface group of genus $1+(g-1)/p=2$,
and $Q$ is a group of automorphisms of a Riemann surface $\mathcal T$ of genus $2$.
Note that the proof is valid so far provided $p>\max\{p_{\lambda},\,
p'_{\lambda},\, s_{\lambda}+1\}$.

\medskip

Now $\Delta/K\cong P\cong C_p$, so $K$ contains the subgroup $\Delta'\Delta^p$
generated by the commutators and $p$-th powers of elements of $\Delta$, and
hence $P$ is isomorphic (as a ${\bf Z}_pQ$-module) to a 1-dimensional quotient
of the 4-dimensional ${\bf Z}_pQ$-module $\Delta/\Delta'\Delta^p \cong
H_1({\mathcal T},{\bf Z}_p) \cong H_1({\mathcal T},{\bf Z})\otimes{\bf Z}_p$. Since $p$
does not divide $s=|Q|$, Maschke's Theorem implies that $H_1({\mathcal T},{\bf
Z}_p)$ is a direct sum of irreducible submodules. Now $H_1({\mathcal T},{\bf
C})=H_1({\mathcal T},{\bf Z})\otimes{\bf C}$ is a direct sum of two $Q$-invariant
subspaces, corresponding under duality to holomorphic and antiholomorphic
differentials in $H^1({\mathcal T},{\bf C})$, and these afford complex conjugate
representations of $Q$ (see [{\bf 13}], for instance). It follows that
$H_1({\mathcal T},{\bf Z}_p)$ is irreducible, or a direct sum of two irreducible
2-dimensional submodules, or a direct sum of four irreducible 1-dimensional
submodules. Since $H_1({\mathcal T},{\bf Z}_p)$ has a $1$-dimensional quotient,
only the last of these three possibilities can arise. We have $p>2$, so by a
theorem of Serre [{\bf 5}, V.3.4] $Q$ is faithfully represented on $H_1({\mathcal
T},{\bf Z}_p)$; thus $Q\leq \GL_1(p)^4\cong(C_{p-1})^4$, so $Q$ is an abelian
group of exponent $e$ dividing $p-1$. Since $\Delta$ is a surface group, the
natural epimorphism $\Gamma\to\Gamma/\Delta\cong Q$ is a surface-kernel
epimorphism.

\medskip

We have $|Q|=s\geq\lambda>6$, and since no Riemann surface of genus $2$ has an
automorphism of order $7$ [{\bf 5}, V.1.11] it follows that $s\geq 8$; thus
$|G|=ps\geq 8(g-1)$ and hence $\sigma\in\Sigma_8$. The elements of $\Sigma_8$
are listed in the Appendix, and by inspection, the only cases in which there is
a surface-kernel epimorphism from $\Gamma$ to an abelian group $Q$ are the
following:
\begin{itemize}
\item[(i)] $\Gamma=\Gamma(2,6,6)$ with $Q\cong C_6\times C_2$,
\vskip4pt
\item[(ii)] $\Gamma=\Gamma(2,5,10)$ with $Q\cong C_{10}$,
\vskip4pt
\item[(iii)] $\Gamma=\Gamma(2,8,8)$ with $Q\cong C_8$ or $C_8\times C_2$.
\end{itemize}
In case~(i) we have $s=12,\, e=6,\, |G|=12p=12(g-1)$ and
$p\equiv 1$ mod~$(3)$, giving conclusion~(i) of Theorem~1(a); in case~(ii) we
have $s=e=10,\, |G|=10p=10(g-1)$ and $p\equiv 1$ mod~$(5)$, giving
conclusion~(ii); in case~(iii) we have $s=e=8$ (so $Q\cong C_8$ since $|Q|=s$),
$|G|=8p=8(g-1)$ and $p\equiv 1$ mod~$(8)$, giving conclusion~(iii). We will
show in Section~5 that in these three cases, $G$ and $\mathcal S$ are as described
in Examples~(i), (ii) and (iii) of Section~3.

\medskip

{\it Case B\/}: $\sigma=(2,3,n)$ for $n\geq 79$. Here it would be sufficient to
argue as in [{\bf 4}, lemma~3.1] to show that $p$ is bounded above, as we will
do in Case~D for $\sigma=(2,5,n)$. However, the bound given by that argument is
rather large, so for future use, when we consider the case $\lambda=8$ in
Section~6, we will provide a more detailed argument which eliminates this case
completely, and which also gives useful information when $n<79$ (Case~A).

\medskip

If $\Gamma=\Gamma(2,3,n)$, then $q=\mu(\Gamma)/4\pi=(n-6)/12n$, and
$|G|=(g-1)/q=12np/(n-6)$. If $p$ does not divide $n-6$ then $n-6$ must divide
$12n$, so $n-6$ divides $12n-12(n-6)=72$ and hence $n\leq 78$, against our
hypothesis. Thus $p$ divides $n-6$, say $n=kp+6$, so $|G|=12np/kp=12n/k$. It is
sufficient to eliminate the case $n=kp+6$ for all $n\geq 79$, but in fact it is
just as easy (and more useful later) to eliminate it on the weaker assumption
that $n\geq 11$.

\medskip

Since there is a surface-kernel epimorphism from $\Gamma(2,3,n)$ to $G$, the
elliptic generator $\gamma_3$ of $\Gamma$ has an image $c$ of order $n$ in $G$,
so $n$ divides $|G|$ and hence $k$ divides $12$. Now $G$ has a cyclic subgroup
$C=\langle c\rangle\cong C_n$ of index $12/k\leq 12$. The kernel of the
transitive action of $G$ on the $12/k$ cosets of $C$ is the core $Z$ of $C$, a
cyclic normal subgroup of $G$, and $c$ induces a permutation $\pi\in S_{12/k}$
of order $l=|C:Z|$ dividing $n$, with at least one fixed point (namely $C$).
Since $c$ centralises $Z$, and the periods $2$ and $3$ are coprime, $Z$ is in
the centre of $G$. By considering the possible cycle-structures of $\pi$, we
show that $n$ is small, thus giving a contradiction.

\medskip

First let $k=1$, so $n=p+6$; thus $n$ is coprime to $6$ since $p\geq 5$, so $G$
is perfect. The cycle-lengths of $\pi$ (in $S_{12}$) are coprime to $6$, so
they must be $1, 5, 7$ or $11$. If there is an $11$-cycle in $\pi$, then $l=11$
and $G/Z$ is a doubly transitive group of degree $12$ and order $12\cdot 11$,
so it is a doubly transitive Frobenius group; however these all have
prime-power degree (since the Frobenius kernel must be elementary abelian [{\bf
7}, XII.9.1]), so there is no $11$-cycle. If there is a $7$-cycle in $\pi$,
then $c$ must fix the remaining $5$ points; however, a little experimentation
shows that elements of order $2$ and $3$ in $S_{12}$ cannot generate a
transitive group and have a $7$-cycle as their product (i.e.~there is no
transitive Hurwitz group of degree $12$), so $\pi$ contains no $7$-cycle. If
there is a $5$-cycle, then all the cycles-lengths are $1$ or $5$, so $l=5$ and
$G/Z\cong\Gamma(2,3,5)\cong A_5\cong \PSL_2(5)$; since $G$ is perfect, $|Z|$
divides the order of the Schur multiplier $|M(G/Z)|=|M(\PSL_2(5))|=2$ [{\bf 6},
V.25.7], so $n$ divides $2l=10$, against our assumption that $n\geq 11$. Hence
there is no 5-cycle, so $\pi=1\in S_{12}$. This means that $C$ is normal in
$G$, but then $G/C$ has order $1$ (since $\gamma_1$ and $\gamma_2$ have coprime
orders) whereas $|G:C|=12$.

\medskip

Now let $k=2$, so $|G:C|=6$, $n=2p+6$ is coprime to $3$, and hence $\pi\in S_6$
has cycle-lengths $1, 2, 4$ or $5$. If there is a $5$-cycle in $\pi$ then $l=5$
and $G/Z$ is a doubly transitive Frobenius group of degree $6$, which is
impossible. If there is a $4$-cycle then $l=4$ and $G/Z\cong\Gamma(2,3,4)\cong
S_4$; now $|G:G'|$ divides $2$ since $n$ is coprime to $3$, and $|S_4:A_4|=2$,
so $Z\leq G'$; hence $|Z|$ divides $|M(S_4)|=2$ [{\bf 6}, V.25.12] and so $n$
divides $2l=8$; however, $n\geq 11$, so there is no $4$-cycle. If there are
$2$-cycles then $l=2$ and $G/Z\cong\Gamma(2,3,2)\cong S_3$, and we obtain a
contradiction again since $Z\leq G'$ and $|M(S_3)|=1$ [{\bf 6}, V.25.12]. Thus
$\pi=1$, so $C$ is normal in $G$ and we obtain a contradiction as in the case
$k=1$.

\medskip

If $k=3$ then $|G:C|=4$, $n=3p+6$ is odd, and hence $\pi\in S_4$ has
cycle-lengths $1$ or $3$. If $\pi$ contains a $3$-cycle, then $l=3$,
$G/Z\cong\Gamma(2,3,3)\cong A_4\cong \PSL_2(3)$ and $Z\leq G'$, so we obtain a
contradiction since $|Z|$ divides $|M(\PSL_2(3))|=2$ [{\bf 6}, V.25.7]. Thus
$\pi=1$ so $C$ is normal and we again have a contradiction.

\medskip

If $k=4$ then $|G:C|=3$ and $n=4p+6$ is coprime to $3$, so $\pi$ has order
$l=1$ or $2$ in $S_3$, giving contradictions as before. If $k=6$ or $12$ then
$C$ is normal in $G$, again giving a contradiction. Thus the case
$\sigma=(2,3,n)$ cannot arise for $n\geq 79$, and we have also shown that when
$p$ divides $n-6$ it cannot arise for $n\geq 11$.

\medskip

{\it Case C\/}: $\sigma=(2,4,n)$ for $n\geq 37$. Here
$q=\mu(\Gamma)/4\pi=(n-4)/8n$, and $|G|=(g-1)/q=8np/(n-4)$. If $p$ does not
divide $n-4$ then $n-4$ must divide $8n$, so $n-4$ divides $8n-8(n-4)=32$ and
hence $n\leq 36$, against our hypothesis. Thus $p$ divides $n-4$, say $n=kp+4$,
so $|G|=8np/kp=8n/k$, where $k$ divides $8$ since the image $c$ of $\gamma_3$
generates a subgroup $C$ of order $n$. As in Case~B, $c$ induces a permutation
$\pi\in S_{8/k}$ on the cosets of $C$, of order $l=|C:Z|$ dividing $n$, with at
least one fixed-point. Here $Z$ is the core of $C$, but unlike in the previous
case, $Z$ is not necessarily central: $G$ induces a group of automorphisms of
$Z$ of order dividing ${\rm hcf}(2,4)=2$.

\medskip

First let $k=1$, so $n=p+4$ is odd and hence the cycle-lengths of $\pi\in S_8$
are $1,3,5$ or $7$. If there is a $7$-cycle, then $G/Z$ is a doubly transitive
Frobenius group of degree $8$; the only such group is $\AGL_1(8)$, and this,
having no elements of order $4$, is not an epimorphic image of $\Gamma(2,4,n)$
for any $n$, so there is no $7$-cycle. If there is a $5$-cycle then the
remaining three points are fixed, so $l=5$; however, by trial and error one can
see that no epimorphic image of $\Gamma(2,4,5)$ can be a transitive group of
degree $8$, so there is no $5$-cycle. By transitivity, $\pi\neq 1$, so $l=3$
and $G/Z$ is an epimorphic image of $\Gamma(2,4,3)\cong S_4$; since $G/Z$ is
transitive of degree $8$ it must be isomorphic to $S_4$, so $G$ is an extension
of $Z\cong C_{n/3}$ by $S_4$. Since $3$ divides $n$ we have $p\equiv 2$
mod~$(3)$, as in conclusion~(iv).

\medskip

If $k=2$ then $n=2p+4$, and $G$ acts transitively on the four cosets of $C$.
Since $\pi$ has a fixed-point, it cannot contain a $4$-cycle. If $\pi$ contains
a $3$-cycle, then $l=3$ and $G/Z$ is an epimorphic image of $\Gamma(2,4,3)\cong
S_4$, whereas there is no such group of order $4l=12$; hence $\pi$ does not
contain a $3$-cycle. By transitivity, $\pi\neq 1$, so $\pi$ is a $2$-cycle,
giving $l=2$ and $G/Z\cong\Gamma(2,4,2)\cong D_4$, so $G$ is an extension of
$Z\cong C_{n/2}$ by $D_4$, as in conclusion~(v).

\medskip

If $k=4$ then $C$ is a normal subgroup of index $2$ in $G$, and
$|G|=2n=8p+8=8g$, as in conclusion~(vi). If $k=8$ then $G=C$ is cyclic, so
$n\leq 4$, against our hypothesis.

\medskip

{\it Case D\/}: $\sigma=(2,5,n)$ for $n\geq 71$. Here
$q=\mu(\Gamma)/4\pi=(3n-10)/20n$, and $|G|=(g-1)/q=20np/(3n-10)$. If $p$ does
not divide $3n-10$ then $3n-10$ must divide $20n$, so $3n-10$ divides $200$ and
hence $n\leq 70$, against our hypothesis. Thus $p$ divides $3n-10$, say
$3n=kp+10$, so $|G|=20np/(3n-10)=20n/k$, where $k$ divides $20$ since $c$
generates a subgroup $C$ of order $n$. As in Case~B, the core $Z$ of $C$ is
central in $G$, since the periods $2$ and $5$ are coprime. Now $C$ has index
$20/k$ in $G$, so $Z$ has index $m$ dividing $(20/k)!$, and since $Z$ is
central, the transfer from $G$ to $Z$ induces the endomorphism $z\mapsto z^m$
of $Z$ [{\bf 6}, IV.2.1]. Since $Z\cap G'$ is in its kernel, and is cyclic, it
has order dividing $m$. Now $|Z:Z\cap G'|=|ZG':G'|$ divides $|G:G'|$, and this
divides $2\cdot 5=10$, so $|Z|$ divides $10m$. Thus $|G|=|G:Z|\cdot |Z|$
divides $10m^2$, so $|G|\leq 10\cdot (20!)^2$ giving $p=g-1=|G|q<3|G|/20\leq
3\cdot(20!)^2/2$.

\medskip

It follows that Theorem~1(a) holds for all $p\geq
c_{\lambda}>\max\{p_{\lambda},\, p'_{\lambda},\, s_{\lambda}+1,\,
3\cdot(20!)^2/2\}$, where $p_{\lambda},\, p'_{\lambda}$ and $s_{\lambda}$ are
as defined in Case~A. Indeed, if $\lambda>20/3$ then only finitely many
signatures $\sigma=(2,5,n)$ satisfy $\mu(\Gamma)\leq 4\pi/\lambda$; we can
enlarge the finite set $\Sigma_{\lambda}$ to include these, so Case~D does not
arise and we can take $c_{\lambda}>\max\{p_{\lambda},\, p'_{\lambda},\,
s_{\lambda}+1\}$. We will use this in Section~6.

\section{\bf Proof of Theorem~1(b) and (c)}

If $G$ and $\mathcal S$ are as in Theorem~1 with $p\geq c_{\lambda}$, then by
comparing the upper bounds $8(g+3)$ and $12(g-1)$ in Theorem~1(a) we see that
$|G|\leq 12(g-1)$ if $g\geq 9$. Applying this to ${\rm Aut}\,({\mathcal S})$
itself, we have $|{\rm Aut}\,({\mathcal S})|\leq 12(g-1)$. Since $|G|$ divides
$|{\rm Aut}\,({\mathcal S})|$ and $|G|\geq\lambda(g-1)>20(g-1)/3$, we therefore
have $G={\rm Aut}\,({\mathcal S})$ and hence $\Gamma=N(K)$. This proves
Theorem~1(b), so we now consider Theorem~1(c).

\medskip

In cases~(i), (ii) and (iii), covered by Case~A of the proof of Theorem~1(a),
we showed that $\Gamma=\Gamma(2,6,6),\, \Gamma(2,5,10)$ or $\Gamma(2,8,8)$, and
$G=\Gamma/K$ is a split extension of $P=\Delta/K\cong C_p$ by an abelian group
$Q=\Gamma/\Delta\cong C_6\times C_2,\, C_{10}$ or $C_8$; to determine $G$ it is
therefore sufficient to find the action of $Q$ by conjugation on $P$. This is
equivalent to the action of $Q$ on a $1$-dimensional quotient of $H_1({\mathcal
T},{\bf Z}_p)$, where ${\mathcal T}={\mathcal H}/\Delta$ has genus $2$, so we now
consider this representation of $Q$.

\medskip

It is easily seen that in each case, $\Gamma$ has a unique normal surface
subgroup $\Delta$ of genus $2$ with abelian quotient: in cases~(i) and (ii) it
is the commutator subgroup $\Gamma'$, and in case~(iii) it is the normal
closure of $\gamma_1\gamma_2^4$, which contains $\Gamma'$ with index $2$. It
follows that ${\mathcal T}$ must be the Riemann surface of genus $2$ described in
Section~3, given by $w^2=z^6-1,\, z^5-1$ and $z^5-z$ respectively, with the
action of $Q$ specified there. The character of $Q$ on $H_1({\mathcal T},{\bf Z})$
is given by $\chi(g)=2-\phi(g)$ for non-identity $g\in Q$, where $g$ fixes
$\phi(g)$ points of $\mathcal T$. By counting fixed-points, and then reducing
$\chi(g)$ mod~$(p)$, one can decompose the character of $Q$ on $H_1({\mathcal
T},{\bf Z}_p)$ into irreducible constituents (see [{\bf 10}] for full details);
in each case there are four distinct $1$-dimensional constituents, which
implies that the four $1$-dimensional submodules of $H_1({\mathcal T},{\bf Z}_p)$
found in the proof of Theorem~1(a) are non-isomorphic. In case~(i), the
automorphism $(w,z)\mapsto (w,e^{\pi i/3}z)$ has eigenvalues $u^{\pm 1}$ and
$u^{\pm 2}$, where $u$ is a primitive $6$-th root of $1$ in ${\bf Z}_p$, while
$(w,z)\mapsto(-w,z)$ has four eigenvalues equal to $-1$; it follows that for
each of the $1$-dimensional quotients, one can find a decomposition $C_6\times
C_2$ of $Q$ so that generators of the two factors have eigenvalues $u$ and $1$.
Thus the action of $Q$ on $P$ is as given in Example~(i) of Section~3, so $G$
is as described there, and hence the four $1$-dimensional quotients correspond
to the four surface subgroups $K_j$ and surfaces ${\mathcal S}_j$ also described
there. The same argument applies in cases~(ii) and (iii), a generator for
$Q\cong C_{10}$ or $C_8$ having the four primitive $10$-th or $8$-th roots of
$1$ as eigenvalues.

\medskip

In cases~(iv) and (v), covered in Case~C of the proof of Theorem~1(a), it was
shown that $\Gamma=\Gamma(2,4,n)$ and $G$ is an extension of $C_{n/3}$ by
$S_4$, or of $C_{n/2}$ by $D_4$. It is shown in [{\bf 9}], as part of a
classification of cyclic coverings of finite rotation groups, that in each case
the surface group $K$ is unique, and that $G$ has the presentation given in
Section~3. It follows that $\mathcal S$ is unique, and is the branched covering of
the sphere described in Section~3.

\medskip

Finally, in case~(vi) it was shown that $\Gamma=\Gamma(2,4,n)$ with
$n=4g=4(p+1)$, and $C=\langle c\rangle$ has index $2$ in $G$. It follows that
the canonical generators $a,b,c$ of $G$ satisfy
$$a^2=b^4=c^n=abc=1,\, c^a=c^u\eqno(5{\cdot}1)$$
where $u^2\equiv 1$ mod~$(n)$. Now $b^{-2}=(ca)^2=c^{1+u}$, so for $b$ to have
order $4$ we require $2(1+u)\equiv 0\not\equiv 1+u$ mod~$(n)$. Hence $n$ is
even, say $n=2m$, and $1+u\equiv m$ mod~$(n)$, so we can take $u=m-1$. Since we
require $u^2\equiv 1$ mod~$(n)$, we need $m$ to be even. Eliminating $b$ from
($5{\cdot}1$) we have the relations
$$a^2=c^n=1,\, c^a=c^{m-1}\eqno(5{\cdot}2)$$
which define the semidihedral group $SD_n$ of order $2n$. Since $|G|=2n$ also,
these are defining relations for $G$, so $G$ has a presentation
$$G=\langle a,b,c\mid a^2=b^4=c^n=abc=1,\, c^a=c^{m-1}\rangle\eqno(5{\cdot}3)$$
where $n=2m$. The uniqueness of this presentation shows that $K$ and hence
$\mathcal S$ are unique.

\section{\bf The case where $\lambda=8$}

The proof of Theorem~1 shows that computing a suitable value of $c_{\lambda}$
is a matter of routine (and tedious) arithmetic: one can use ($4{\cdot}2$) to
determine the finite set $\Sigma_{\lambda}$, then find the values of
$p_{\lambda},\, p'_{\lambda}$ and $s_{\lambda}$ by inspection, and take
$$c_{\lambda}>\max\Bigl\{p_{\lambda},\, p'_{\lambda},\, s_{\lambda}+1,\,
{3\over2}(20!)^2\Bigr\}.$$ If required, one can use additional, more
sophisticated arguments to provide a lower value of $c_{\lambda}$ by
strengthening the arguments in the proof of Theorem~1 so that they apply to
smaller primes. (Here it is sufficient to restrict attention to Theorem~1(a),
since the proofs of parts~(b) and (c) follow as in Section~5.) To illustrate
this, we will show that when $\lambda=8$ the proof of Theorem~1(a) gives a
lower bound $p>85$, and then we will use additional arguments to extend this to
$p\geq 17$. (The bound cannot be reduced further: when $p=13$ we find a Hurwitz
group $\PSL_2(13)$ of genus $g=14$, a quotient of $\Gamma(2,3,7)$ of order
$84(g-1)$.)

\medskip

Let $\lambda=8$, so $|G|\geq 8(g-1)$ and $\mu(\Gamma)\leq\pi/2$. By detailed
analysis of ($4{\cdot}2$) we obtain the following three types of signatures
$\sigma$ (see [{\bf 2}] for $\mu(\Gamma)<\pi/2$):

\begin{itemize}
\item[(I)] $(2,5,n)$ for $5\leq n\leq 20$; $(2,6,n)$ for $6\leq n\leq 12$;
$(2,7,n)$ for $7\leq n\leq 9$; $(2,8,8)$; $(3,3,n)$ for $4\leq n\leq 12$;
$(3,4,n)$ for $4\leq n\leq 6$; $(2,2,2,n)$ for $n=3,4$;

\vskip2pt

\item[(II)] $(2,3,n)$ for $n\geq 7$;

\vskip2pt

\item[(III)] $(2,4,n)$ for $n\geq 5$.
\end{itemize}

\vskip2pt

\noindent Then $\Sigma_8$ is the finite set consisting of all $41$ signatures
of type I, together with $(2,3,n)$ for $7\leq n\leq 78$, and $(2,4,n)$ for
$5\leq n\leq 36$. The values of $q$ for $\sigma\in\Sigma_8$ are given in the
Table of Signatures in the Appendix. By inspection, the largest prime dividing
any $r$ for $\sigma\in\Sigma_8$ is $71$, which occurs for $\sigma=(2,3,77)$
with $q=71/924$, so we need $p>p_{\lambda}=71$. The Table shows that for those
$\sigma$ with $r=1$, the maximum value of $s$ is $84$, arising from
$\sigma=(2,3,7)$ with $q=1/84$, so taking $s_{\lambda}=84$ we require
$p>s_{\lambda}+1=85$. Among the signatures with $r=1$, the largest prime
dividing any of the elliptic periods $m_i$ is $13$, arising from
$\sigma=(2,3,78)$ with $q=1/13$, so taking $p'_{\lambda}=13$ we require
$p>p'_{\lambda}=13$. It follows that when $\lambda=8$ the proof of Theorem~1(a)
is valid for all $p>85$.

\medskip

We will now use special arguments to extend this proof to the remaining primes
$p\geq 17$, namely $p=17,\, 19,\, 23,\, 29,\, 31,\, 37,\, 41,\, 43,\, 47,\,
53,\, 59,\, 61,\, 67,\, 71,\, 73,\, 79$ and $83$. Recall from Section~4 that
$|G|=ps/r$ with $r$ and $s$ coprime, so $r=1$ or $p$.

\medskip

\noindent(i) First let $p=83$ or $79$. The Table shows that there are no
signatures $\sigma\in\Sigma_8$ with $r=83$ or $79$, so $r=1$ and $|G|=ps$. The
entries in the Table with $r=1$ all have $s\leq 48$ or $s=84$, so Sylow's
Theorems imply that either $G$ has a normal Sylow $p$-subgroup $P$ of order
$p$, or $p=83$ and $G$ has $s=84$ Sylow $p$-subgroups. In the first case, the
proof continues as in Section~4, while the second case is eliminated by the
following result.

\medskip

L{\small EMMA} 4. {\it Let $p$ be a prime and let $G$ be a group which does not
have a normal Sylow $p$-subgroup. If $|G|=p(p+1)$ then $p=2$ or $p$ is a
Mersenne prime, and if $|G|=2p(p+1)$ then $p\leq 5$ or $p$ is a Mersenne
prime.}

\medskip

{\it Proof.} By Sylow's Theorems, the number $n_p$ of Sylow $p$-subgroups of
$G$ divides $|G|$ and satisfies $n_p\equiv 1$ mod~$(p)$, so $n_p=1$ or $p+1$.
By our hypothesis $n_p>1$, so $n_p=p+1$. By Sylow's Theorems, $G$ acts on its
$p+1$ Sylow $p$-subgroups by conjugation as a transitive permutation group
$\widetilde G$. In this action, a Sylow $p$-subgroup $P$ normalises itself, but
normalises no other Sylow $p$-subgroup $P^*$ (for otherwise $PP^*$ would be a
subgroup of $G$ of order $p^2$, contradicting Lagrange's Theorem). A generator
of $P$ therefore fixes $P$ and has a single cycle of length $p$ on the
remaining Sylow $p$-subgroups, so $G$ acts as a doubly transitive group. If
$|G|=p(p+1)$, the stabiliser of two points is trivial, so $\widetilde G$ is
sharply $2$-transitive, that is, a doubly transitive Frobenius group. Such
groups all have prime-power degree [{\bf 7}, XII.9.1], so $p=2$ or $2^e-1$ for
some $e$. If $|G|=2p(p+1)$ then either $|\widetilde G|=p(p+1)$, and the
preceding argument applies, or $|\widetilde G|=2p(p+1)$ with a point-stabiliser
acting on the remaining $p$ points as $D_p$ (the only transitive group of
degree $p$ and order $2p$); in this case, $\widetilde G$ is a Zassenhaus group
with two-point stabilisers of even order ($=2$), and Zassenhaus showed that
such groups have two-point stabilisers of order at least $(p-2)/2$ [{\bf 7},
XI.1.10], so $p\leq 5$.

\medskip

\noindent(ii) Now suppose that $p=73,\, 71,\, 67,\, 61,\, 59,\, 53,\, 43$ or
$37$. The case $r=1$ is dealt with as in (i): the Table shows that $s=84,\,
48,\, 40$ or $s\leq 36$, none with a factor $kp+1>1$, so $G$ has a normal Sylow
$p$-subgroup. If $r=p$, the Table shows that $\sigma=(2,3,n)$ with
$s/r=12n/(n-6)$, so $p$ divides $n-6$. However, in the proof of Case~B in
Section~4, this possibility was eliminated for all $p\geq 5$.

\medskip

\noindent(iii) If $p=47$ then we can deal with the case $r=1$ as we did in (i)
for $p=84$, since $47$ is not a Mersenne prime. Hence we can assume that $r=p$,
so the Table gives $\sigma=(2,3,53)$ or $(2,5,19)$. The signature $(2,3,53)$ is
eliminated as in (ii), since $p$ divides $n-6$. If $\sigma=(2,5,19)$ then
$|G|=2^2\cdot 5\cdot 19$, so $G$ has a normal Sylow $19$-subgroup (by Lemma~4),
with a solvable quotient by Burnside's $p^aq^b$ Theorem [{\bf 6}, V.7.3], so
$G$ is solvable; however, the periods of $\Gamma$ are mutually coprime, so $G$
is perfect, giving a contradiction. Similar arguments apply to $p=41$. If $r=1$
then $G$ has a normal Sylow $41$-subgroup, as required, since the only $s$ in
the Table divisible by some $kp+1>1$ is $s=84$, with $\sigma=(2,3,7)$ and
$|G|=84\cdot 41$, eliminated by Lemma~4. If $r=p$ then either
$\sigma=(2,3,47)$, or $\sigma=(2,5,17)$ and $|G|=2^2\cdot 5\cdot 17$, both
eliminated as for $p=47$.

\medskip

\noindent(iv) Now let $p=31$. No $kp+1>1$ divides any $s$ in the Table with
$r=1$, so we may assume that $r=p$. The Table then gives $\sigma=(2,3,37),\,
(2,3,68),\, (2,4,35)$ or $(2,7,9)$. The first two are eliminated as in (ii)
since $p$ divides $n-6$. If $\sigma=(2,7,9)$ then $|G|=2^2\cdot 3^2\cdot 7$, so
$G$ has a normal Sylow $7$-subgroup and is therefore solvable, contradicting
the fact that $\Gamma$ is perfect. In the remaining case $\sigma=(2,4,35)$,
with $|G|=2^3\cdot 5\cdot 7=280$, note that $|G^{\rm ab}|$ divides
$|\Gamma^{\rm ab}|=2$. Let $n_2, n_5$ and $n_7$ denote the numbers of Sylow
$2$-, $5$- and $7$-subgroups of $G$, so $n_5=1$ or $56$, and $n_7=1$ or $8$. If
$n_5=56$ and $n_7=8$ then there are $1+56\cdot 4+8\cdot 6=273$ elements of
Sylow $5$- or $7$-subgroups, so $n_2=1$; thus there is a normal Sylow
$2$-subgroup $N$, so $G/N$, having order $5\cdot 7$, must be abelian,
contradicting the fact that $|G^{\rm ab}|$ divides $2$. Hence $n_5=1$ or
$n_7=1$. If $n_7=1$ there is a normal Sylow $7$-subgroup $N$, and $G/N$ (of
order $2^3\cdot 5$) has a normal Sylow $5$-subgroup, so $G$ has a normal
subgroup $M$ of order $5\cdot 7$; thus $\Gamma$ maps onto a group $G/M$ of
order $2^3$, which is clearly false. Hence $n_5=1$, so $G$ has a normal Sylow
$5$-subgroup $N$; by the previous argument, $G/N$ cannot have a normal Sylow
$7$-subgroup, so it has eight of them and hence has a normal Sylow
$2$-subgroup; thus $G$ has a normal subgroup of index $7$, so $\Gamma$ maps
onto $C_7$, which is false.

\medskip

\noindent(v)  If $p=29$ then we can deal with the case $r=1$ as we did in (i)
for $p=84$, since $p$ is not a Mersenne prime. We therefore have $r=p$, so
$\sigma=(2,5,13)$ with $|G|=2^2\cdot 5\cdot 13$, or $(2,3,35)$ with
$|G|=2^2\cdot 3\cdot 5\cdot 7$, or $(2,3,64)$ with $|G|=2^7\cdot 3$, or
$(2,4,33)$ with $|G|=2^3\cdot 3\cdot 11$. In the first case, $G$ has a normal
Sylow $13$-subgroup and is solvable, contradicting the fact that $\Gamma$ is
perfect. The second and third cases are eliminated as in Case~B of the proof of
Theorem~1(a), since $p$ divides $n-6$. The last case is dealt with in Case~C,
where $\sigma=(2,4,n)$ and $n=kp+4$ with $k=1$: we showed there that $G$ is an
extension of $C_{n/3}=C_{11}$ by $S_4$, as in conclusion~(iv) of Theorem~1(a).

\medskip

\noindent(vi) Let $p=23$. If $r=1$, the only values of $s$ with factors
$kp+1>1$ are $24$ and $48$, both eliminated by Lemma~4, so $n_p=1$ as required.
If $r=p$ then the Table gives $\sigma=(2,3,29),\, (2,3,52),\, (2,3,75),\,
(2,4,27)$ or $(2,5,11)$. The first three are eliminated as in Case~B, since $p$
divides $n-6$. The fourth case is dealt with in Case~C, with $n=p+4$, where we
showed that $G$ is an extension of $C_9$ by $S_4$, as in conclusion~(iv) of
Theorem~1. In the final case, $|G|=2^2\cdot 5\cdot 11$, so $G$ has a normal
Sylow $11$-subgroup and is solvable, contradicting the fact that $\Gamma$ is
perfect.

\medskip

\noindent(vii) Let $p=19$. If $r=1$, the only values of $s$ with factors
$kp+1>1$ are $20$ and $40$, both eliminated by Lemma~4, so $n_p=1$ as required.
If $r=p$ then the Table gives $\sigma=(2,3,25),\, (2,3,44),\, (2,3,63),\,
(2,4,23)$ or $(2,5,16)$. The first three are eliminated as in Case~B, since $p$
divides $n-6$, and the fourth case is eliminated in Case~C, with $n=p+4$, since
$19\not\equiv 2$ mod~$(3)$. In the final case, $|G|=160=2^5\cdot 5$, so $G$ is
solvable. Now $\Gamma/\Gamma'\cong C_2$, with $\Gamma'\cong\Gamma(5,5,8)$, and
so $\Gamma'/\Gamma''\cong C_5$, giving $G/G'\cong C_2$ and $G'/G''\cong C_5$,
and hence $|G''|=2^4$. By Maschke's Theorem, the $G'/G''$-module
$G''/\Phi(G'')$ is a sum of irreducible submodules (where $\Phi$ denotes the
Frattini subgroup). Since $G''$ has no quotient of order $2$, there can be no
$1$-dimensional summands, and $C_5$ has no $2$- or $3$-dimensional irreducible
modules over ${\bf Z}_2$, so $G''/\Phi(G'')$ is $4$-dimensional and has order
$2^4$. Thus $\Phi(G'')=1$, so $G''$ is elementary abelian. This means that a
Sylow $2$-subgroup of $G$ has exponent dividing $4$, contradicting the fact
that $G$ contains an element of order $16$.

\medskip

\noindent(viii) Let $p=17$. If $r=1$, the only values of $s$ with factors
$kp+1>1$ are $18$ and $36$, both eliminated by Lemma~4, so $n_p=1$ as required.
If $r=p$ then the Table gives $\sigma=(2,3,23),\, (2,3,40),\, (2,3,57),\,
(2,3,74),\, (2,4,21)$ or $(2,5,9)$. The first four are eliminated as in Case~B,
since $p$ divides $n-6$, while $(2,4,21)$ is dealt with in Case~C, where
$n=p+4$ and $G$ is an extension of $C_7$ by $S_4$, as in conclusion~(iv) of
Theorem~1. If $\sigma=(2,5,9)$ then $|G|=180$; now $\Gamma$ is perfect, whereas
there is no perfect group of order $180$ (since $|M(A_5)|=2$, the only
extension of $C_3$ by $A_5$ is the direct product).

\section{\bf Proof of Theorem 2}

We can now complete the proof of Theorem~2, which we began in Corollary~3 in
Section~3.

\medskip

\noindent(a) If $p\equiv 1$ mod~$(3)$ then Example~(i) of Section~3 gives
$N(g)\geq 12(g-1)$ for all $g=p+1$, and the case $\lambda=8$ of Theorem~1,
considered in Section~6, shows that we have equality if $p\geq 17$.

\medskip

\noindent(b) If $p\equiv 11$ mod~$(15)$ then $p\equiv 1$ mod~$(5)$, so
Example~(ii) gives $N(g)\geq 10(g-1)$; since $p\not\equiv 1$ mod~$(3)$, we
again have equality for $p\geq 17$.

\medskip

\noindent(c) If $p\equiv 2, 8$ or $14$ mod~$(15)$ then $p\equiv 2$ mod~$(3)$,
so Example~(iv) gives $N(g)\geq 8(g+3)$; since $p\not\equiv 1$ mod~$(3)$ and
$p\not\equiv 1$ mod~$(5)$, we have equality for $p\geq 17$.

\medskip C{\small OROLLARY} 5. {\it If $g=p+1$ for some prime $p$ then
$N(g)\geq 8(g+3)$, and this bound is attained for all $p\geq 17$ such that
$p\equiv 2, 8$ or $14$ {\rm mod}~$(15)$.}

\medskip

It follows from Dirichlet's Theorem that this bound is attained infinitely many
times.

\bigskip\bigskip

\centerline{REFERENCES}

\bigskip

\begin{list}{}
{
\setlength{\itemsep}{4pt}
\setlength{\leftmargin}{0pt}
\setlength{\itemindent}{20pt}
}
\item[{\bf 1.}] R.~D.~A{\small CCOLA}. On the number of automorphisms of a
closed Riemann surface. {\it Trans. Amer.~Math.~Soc.} {\bf 131} (1968),
398--408.
\item[{\bf 2.}] R.~D.~A{\small CCOLA}. {\it Topics in the Theory of Riemann
Surfaces}, Lecture Notes in Math. {\bf 1595} (Springer-Verlag, 1994).

\item[{\bf 3.}] M.~B{\small ELOLIPETSKY} and G.~A.~J{\small ONES}. A bound
for the number of automorphisms of an arithmetic Riemann surface. {\it
Math.~Proc.~Cambridge Phil.~Soc.}, to appear.

\item[{\bf 4.}] M.~D.~E.~C{\small ONDER} and R.~K{\small ULKARNI}.
Infinite families of automorphism groups of Riemann surfaces. In {\it Groups
and Geometry\/} (eds.~C.~Maclachlan and W.~J.~Harvey), pp.~47--56. {\it
Lond.~Math.~Soc.~Lect.~Note Ser.} {\bf 173} (Cambridge University Press, 1992).

\item[{\bf 5.}] H.~M.~F{\small ARKAS} and I.~K{\small RA}. {\it Riemann
Surfaces}. Grad.~Texts in Math. {\bf 71} (Springer-Verlag, 1980).

\item[{\bf 6.}] B.~H{\small UPPERT}. {\it Endliche Gruppen I\/}
(Springer-Verlag, 1967).

\item[{\bf 7.}] B.~H{\small UPPERT} and N.~B{\small LACKBURN}. {\it Finite
Groups III\/} (Springer-Verlag, 1982).

\item[{\bf 8.}] G.~A.~J{\small ONES} and D.~S{\small INGERMAN}. {\it
Complex Function Theory: an Algebraic and Geometric Viewpoint\/} (Cambridge
University Press, 1987).

\item[{\bf 9.}] G.~A.~J{\small ONES} and D.~B.~S{\small UROWSKI}. Regular
cyclic coverings of the Platonic maps. {\it European J.~Combinatorics\/} {\bf
21} (2000), 333--345.

\item[{\bf 10.}] M.~K{\small AZAZ}. Finite Groups and Surface Coverings.
Ph.D.~thesis. University of South\-ampton (1997).

\item[{\bf 11.}] C.~M{\small ACLACHLAN}. A bound for the number of
automorphisms of a compact Riemann surface. {\it J.~London Math.~Soc.} {\bf 44}
(1968), 265--272.

\item[{\bf 12.}] G.~R{\small IERA} and R.~R{\small ODR\'IGUEZ}. Riemann
surfaces and abelian varieties with an automorphism of prime order. {\it Duke
Math.~J.} {\bf 69} (1993), 199--217.

\item[{\bf 13.}] C.-H.~S{\small AH}. Groups related to compact Riemann
surfaces. {\it Acta Math.} {\bf 123} (1969), 13--42.

\item[{\bf 14.}] D.~S{\small INGERMAN}. Finitely maximal Fuchsian groups.
{\it J.~London Math.~Soc.}~(2) {\bf 6} (1972), 29--38.

\item[{\bf 15.}] M.~S{\small TREIT} and J.~W{\small OLFART}. Characters and
Galois invariants of regular dessins. {\it Revista Matematica Complutense\/}
{\bf 13} (2000), 1-33.

\end{list}

%\vfil\eject
\bigskip\bigskip

\vbox{ \centerline{Appendix. {\it Table of signatures $\sigma\in\Sigma_8$, with
values of $s/r$}}

\bigskip

\centerline{ \vbox{\tabskip=0pt \offinterlineskip
\def\tablerule{\noalign{\hrule}}
\halign to111pt{\strut#& \vrule#\tabskip=1em plus2em&
  #& \vrule#& \hfil#\hfil& \vrule #&
  #\hfil& \vrule#\tabskip=0pt\cr\tablerule
&&\omit\hidewidth $\sigma$\hidewidth&& \omit\hidewidth
$s/r$\hidewidth&\cr\tablerule &&$(2,5,5)$&& 20      &\cr\tablerule
&&$(2,5,6)$&& 15      &\cr\tablerule &&$(2,5,7)$&& 140/11  &\cr\tablerule
&&$(2,5,8)$&& 80/7    &\cr\tablerule &&$(2,5,9)$&& 180/17  &\cr\tablerule
&&$(2,5,10)$&& 10     &\cr\tablerule &&$(2,5,11)$&& 220/23 &\cr\tablerule
&&$(2,5,12)$&& 120/13 &\cr\tablerule &&$(2,5,13)$&& 260/29 &\cr\tablerule
&&$(2,5,14)$&& 35/4   &\cr\tablerule &&$(2,5,15)$&& 60/7   &\cr\tablerule
&&$(2,5,16)$&& 160/19 &\cr\tablerule &&$(2,5,17)$&& 340/41 &\cr\tablerule
&&$(2,5,18)$&& 90/11  &\cr\tablerule &&$(2,5,19)$&& 380/47 &\cr\tablerule
&&$(2,5,20)$&& 8      &\cr\tablerule &&$(2,6,6)$&& 12      &\cr\tablerule
&&$(2,6,7)$&& 21/2    &\cr\tablerule &&$(2,6,8)$&& 48/5    &\cr\tablerule
&&$(2,6,9)$&& 9       &\cr\tablerule &&$(2,6,10)$&& 60/7   &\cr\tablerule
&&$(2,6,11)$&& 33/4  &\cr\tablerule \noalign{\smallskip}\cr}} \quad\quad
\vbox{\tabskip=0pt \offinterlineskip
\def\tablerule{\noalign{\hrule}}
\halign to187pt{\strut#& \vrule#\tabskip=1em plus2em&
  #& \vrule#& \hfil#\hfil& \vrule#&
  #\hfil& \vrule#\tabskip=0pt\cr\tablerule
&&\omit\hidewidth $\sigma$\hidewidth&& \omit\hidewidth
$s/r$\hidewidth&\cr\tablerule &&$(2,6,12)$&& 8     &\cr\tablerule &&$(2,7,7)$&&
28/3   &\cr\tablerule &&$(2,7,8)$&& 112/13 &\cr\tablerule &&$(2,7,9)$&& 252/31
&\cr\tablerule &&$(2,8,8)$&& 8      &\cr\tablerule &&$(3,3,4)$&& 24
&\cr\tablerule &&$(3,3,5)$&& 15     &\cr\tablerule &&$(3,3,6)$&& 12
&\cr\tablerule &&$(3,3,7)$&& 21/2   &\cr\tablerule &&$(3,3,8)$&& 48/5
&\cr\tablerule &&$(3,3,9)$&& 9      &\cr\tablerule &&$(3,3,10)$&& 60/7
&\cr\tablerule &&$(3,3,11)$&& 33/4  &\cr\tablerule &&$(3,3,12)$&& 8
&\cr\tablerule &&$(3,4,4)$&& 12     &\cr\tablerule &&$(3,4,5)$&& 120/13
&\cr\tablerule &&$(3,4,6)$&& 8      &\cr\tablerule &&$(2,2,2,3)$&& 12
&\cr\tablerule &&$(2,2,2,4)$&& 8    &\cr\tablerule &&$(2,3,n),\; 7\leq n\leq
78$ && $12n/(n-6)$ &\cr\tablerule &&$(2,4,n),\; 5\leq n\leq 36$ && $8n/(n-4)$
&\cr\tablerule \cr \noalign{\smallskip}\cr}} } }

\bigskip\bigskip

\end{document}